\documentclass[letterpaper, 10 pt, conference]{ieeeconf}  

\IEEEoverridecommandlockouts                              

\usepackage{amsmath}																									%
\usepackage{amsfonts}																									%
\usepackage{amssymb}																									%
\usepackage{algorithmicx,algpseudocode}															%
\usepackage{algorithm}																						%
\usepackage{tcolorbox} 																						
\usepackage{nicefrac}   																			
\usepackage[]{units}																				
\usepackage[latin1]{inputenc}																						%
\usepackage{type1cm}																									%
\usepackage{type1ec}																										%
\usepackage[T1]{fontenc}																								%
\usepackage{dsfont}																										%
\usepackage{booktabs}																									%
\usepackage{array}																											%
\usepackage{enumerate}																									%
\usepackage{makeidx}																										%
\usepackage{fancyhdr}																									%
\usepackage{bm}																												%
\usepackage{lastpage}																									%
\usepackage{multicol}																										%
\usepackage{amsmath, amssymb, amsfonts}																	%
\usepackage{ifthen}																										%
\usepackage{ifpdf}																											%
\usepackage{framed}																										%
\usepackage[amsmath,thmmarks,framed,hyperref]{ntheorem}										%
\DeclareMathAlphabet{\mathpzc}{OT1}{pzc}{m}{it}														%
\usepackage{amsmath,amsfonts,amssymb,booktabs}													%
\usepackage{amstext}																										%
\usepackage{psfrag} 																										%
\usepackage{float} 																											%
\usepackage{cite} 																											%
\usepackage{subfloat} 																									%
\usepackage{stackrel} 																									%
\usepackage{graphicx} 																									%
\usepackage{textcomp} 																									%
\usepackage{color} 																											%
\usepackage{bbm} 																											%
\usepackage{lettrine}																										%
\usepackage{dsfont} 																										%
\usepackage{multirow}																									%
\usepackage[section,verbose]{placeins} 																		%
\usepackage{url}																												%
\usepackage{tikz}																											%
\usetikzlibrary{matrix}																										%
\usetikzlibrary{arrows,automata,matrix,fit,positioning,calc}             								%
\usepackage{pgfplots}																									%
\usepackage[latin1]{inputenc}																						%
\usepackage{placeins}																										%
\usepackage{cite}																											%
\usepackage{pgfplots}
\usetikzlibrary{arrows,automata,matrix,fit,positioning,calc}
\usepackage{verbatim}

\usepackage{mathtools}
\usepackage{pst-node}
\usepackage{lipsum}

\usepackage{balance}

\newcounter{thmCounter}
\newtheorem{lemma}[thmCounter]{\bfseries Lemma}

\usepackage{lipsum}

\usepackage{xcolor}

\usepackage{graphicx}

\title{\LARGE \bf
Sparsity-Promoting Iterative Learning Control for Resource-Constrained Control Systems
}

\author{Burak Demirel, Euhanna Ghadimi, and Daniel E. Quevedo
\thanks{B. Demirel and D. E. Quevedo are with the Chair for Automatic Control (EIM-E), Paderborn University, Warburger Stra\ss e 100, 33098, Paderborn, Germany
        {\tt\small burak.demirel@protonmail.com}
        {\tt\small dquevedo@ieee.org}.
       E. Ghadimi is with the Scania Group, SE-15187 S{\"o}dert{\"a}lje, Sweden, 
		{\tt\small euhanna.ghadimi@scania.com}.
        }
}

\begin{document}

\maketitle
\thispagestyle{empty}
\pagestyle{empty}

\begin{abstract}
We propose novel iterative learning control algorithms to track a reference trajectory in resource-constrained control systems. In many applications, there are constraints on the number of control actions, delivered to the actuator from the controller, due to the limited bandwidth of communication channels or battery-operated sensors and actuators. We devise iterative learning techniques that create sparse control sequences with reduced communication and actuation instances while providing sensible reference tracking precision. Numerical simulations are provided to demonstrate the effectiveness of the proposed control method. 
\end{abstract}


\section{Introduction}\label{sec: intro}
A multitude of techniques are now available in the literature for precise control of mechatronic systems; see, e.g.,~\cite{GGD+:17}. Iterative Learning Control (ILC) is one of the well-known techniques for accurately tracking reference trajectories in industrial systems, which repetitively executes a predefined operation over a finite duration; see, e.g.,~\cite{Moo:93, BTA:06, Owe:16}. The key idea of iterative learning control relies on the use of the information gained from previous trails to update control inputs to be applied to the plant on the next trial. Iterative learning control was first introduced by Arimoto et al.~\cite{AKM:84} to achieve high accuracy control of mechatronic systems. Since the original work was published in 1984, it has been successfully practiced in various areas, including additive manufacturing machines~\cite{BHA+:11}, robotic arms~\cite{VDS:11}, printing systems~\cite{BOK+:14}, electron microscopes~\cite{CTL+:09}, and wafer stages~\cite{MCT:07}.

Modern industrial systems, which employ a large number of spatially distributed sensors and actuators to monitor and control physical processes, suffer from resource -- control, communication, and computation -- constraints. To provide a guaranteed performance or even preserve the stability of the closed-loop systems, it is necessary to take these limitations into account while designing and implementing control algorithms. Sometimes the limited bandwidth of legacy  communication networks  imposes a constraint on the rate of data transmissions. Besides, when the feedback loop is closed over wireless networks, a further resource constraint becomes apparent due to the use of battery-powered sensors and actuators~\cite{Dem:15}. The reduced actuator activity also prolongs the lifetime of actuators or  improves the fuel efficiency. Therefore, it is desirable to have either sparse or sporadically changing control commands to reduce the use of actuators.

Sparsity-promoting techniques, which is borrowed from compressive sensing literature, have been successfully applied to a number of control problems to tackle the resource constraints mentioned above; see e.g.,~\cite{GaM:12, HGM:13, CGH+:13, NQO:14, NOQ:16, NQN:16}. The authors of~\cite{GaM:12, HGM:13, CGH+:13} modified the original model predictive control cost with an $\ell_{1}^{}$-penalty term to promote the sparsity in the control input trajectory. The authors of~\cite{NQO:14, NOQ:16, NQN:16} designed energy-aware control algorithms to limit the actuator activity while providing an attainable control performance. Their design is also based on sparse optimization using $\ell_{1}^{}$-norm. To the best of our knowledge, the design of iterative learning control algorithms for resource-constrained systems has not been addressed in the literature and is  subject of this paper.

\textit{Contributions.}
In this paper, we develop a Sparsity-promoting Iterative Learning Control (S-ILC) technique for resource-constrained control systems. The main departure from the standard ILC approach is that we introduce a  regularization term into the usual $\ell_2$-norm cost functions to render the resulting control inputs sparse. The sparsity here is in the cardinality of changes in control values applied to a finite horizon.  Moreover, we include additional constraints to model the practical limits on the magnitude of applied control signals. The resulting control problem is then solved using a backward-forward splitting method which trades off between minimizing the tracking error and finding a sparse control input that optimizes the cost with respect to regularizer term. We demonstrate the monotonic convergence of the technique in lack of modeling imperfections. Moreover, we develop an accelerated algorithm to reduce the number of  trials required for S-ILC to converge to optimality. 

\textit{Outline.}
The remainder of this paper is organized as follows: Section~\ref{sec:problem_formulation} introduces the problem definition. Section~\ref{sec:SILC} presents the sparse iterative learning control problem and associated algorithms to solve it. A numerical study is performed in Section~\ref{sec:numerical}. Finally, Section~\ref{sec:conclusion} presents concluding remarks. The appendix provides  proofs of the main results

\textit{Notation.}
The $n$-dimensional real space is represented by $\mathbb{R}^n$. $\mathbb{E}$ denotes a finite dimensional euclidean space with inner product $\langle \cdot, \cdot \rangle$. For $u\in \mathbb{R}^n$, its $\ell_1$ and $\ell_2$ norms are
\begin{align*}
	\parallel u \parallel_{1}^{} :=  \sum_{i=1}^{n} \vert u_{i}^{} \vert \;, \quad
 	\parallel u \parallel_{}^{} :=  \Bigg( \sum_{i=1}^{n} u_{i}^{2} \Bigg)_{}^{\frac{1}{2}}.
\end{align*}
The spectral radius of the real square matrix $M\in \mathbb{R}^{n\times n}$ is denoted by $\rho(M)$. The Euclidean projection of $u\in \mathbb{R}^n$ into the compact convex set $\mathcal{U}$ is denoted by $\Pi_\mathcal{U}^{}(\cdot)$.

\section{Problem Formulation}\label{sec:problem_formulation}

\subsection{System model}
We consider the following discrete-time, single-input-single-output (SISO), stable, linear time-invariant (LTI) system $P(z)$ with state space representation:
\begin{align}
	x_{k}^{}[t+1] =&\; A x_{k}^{}[t] + B u_{k}^{}[t] \;, \label{eqn:system_model_1a} \\
	y_{k}^{}[t] =&\; C x_{k}^{}[t] \;, \label{eqn:system_model_1b}
\end{align}
where $t\in\mathbb{N}_{0}^{}$ is the time index (i.e., sample number), $k\in\mathbb{N}_{0}^{}$ is the iteration number, $x_{k}^{}[t]\in\mathbb{R}_{}^{n}$ is the state variable, $u_{k}^{}[t]\in\mathbb{R}$ is the control input, $y_{k}^{}[t]\in\mathbb{R}$ is the output variable, and $A$, $B$ and $C$ are matrices of appropriate dimensions. The initial condition $x[0] = x_{0}^{}$ is also assumed to be given, and these initial conditions are the same at the beginning of each trial. The input-output behavior of the system  in~\eqref{eqn:system_model_1a} and~\eqref{eqn:system_model_1b}, can be described via a convolution of the input with the impulse response of the system:
\begin{align}
	y_{k}^{}[t] = CA_{}^{t}x_{0}^{} + \sum_{\tau=0}^{t-1} CA_{}^{t-\tau-1}Bu_{k}^{}[\tau] \;. 
	\label{eqn:impulse_response}
\end{align}
The coefficients $CA_{}^{t}B$ for any $t \in\{0, 1,\cdots,\mathrm{T}\}$ are referred to as the Markov parameters of the plant $P(z)$, provided in~\eqref{eqn:system_model_1a} and~\eqref{eqn:system_model_1b}.

\subsection{Lifted system model}
Since we focus on a finite trial length $\mathrm{T}$, it is possible to evaluate~\eqref{eqn:impulse_response} for all $t\in\{0,1,\cdots,\mathrm{T}\}$ and, similar to~\cite{OHD:09}, write its lifted version as
\begin{align}
	y_{k}^{} = G u_{k}^{} + d \;,
\end{align}
where
\begin{align*}
	G =&\;
	\begin{bmatrix}
		CA_{}^{t_{}^{*}-1}B 	& 0 		& \cdots & 0 \\
		CA_{}^{t_{}^{*}}B 	& CA_{}^{t_{}^{*}-1}B 	& \cdots & 0 \\
		\vdots & \vdots  & \ddots & \vdots \\
		CA_{}^{\mathrm{T}-1}B & CA_{}^{\mathrm{T}-2}B & \cdots & CA_{}^{t_{}^{*}-1}B
	\end{bmatrix} \;, \\
	d =&\; 
	\begin{bmatrix}
		CA_{}^{t_{}^{*}} x_{0}^{} & CA_{}^{t_{}^{*}+1} x_{0}^{} & \cdots & CA_{}^{\mathrm{T}} x_{0}^{}
	\end{bmatrix}_{}^{\top} \;.
\end{align*}
The vectors of inputs and output series are defined as 
\begin{align*}
	u_{k}^{} =&\; 
	\begin{bmatrix}
		u_{k}^{}[0] & u_{k}^{}[1] & \cdots & u_{k}^{}[\mathrm{T}-t_{}^{*}] 
	\end{bmatrix}_{}^{\top} \;, \\
	y_{k}^{} =&\; 
	\begin{bmatrix}
		y_{k}^{}[t_{}^{*}] & y_{k}^{}[t_{}^{*}+1] & \cdots & y_{k}^{}[\mathrm{T}] 
	\end{bmatrix}_{}^{\top} \;.
\end{align*}
The relative degree of the transfer function $P(z)$ is denoted by $t_{}^{*} > 0$. Notice that the matrix $G$ has a Toeplitz structure.

\subsection{Trajectory tracking problem}
In this paper, we focus on the reference trajectory tracking problem. It is assumed that a reference trajectory $r[t]$ is given over a finite time-interval between $0$ and $\mathrm{T}$. The objective is, here, to determine a control input trajectory $\{u[t]\}_{t=0}^{\mathrm{T}-t_{}^{*}}$ that minimizes the tracking error:
\begin{align}
	\parallel e \parallel_{}^{2} \; \triangleq \; \parallel r - y \parallel_{}^{2} \; = \; \parallel r - G u \parallel_{}^{2} \;,
	\label{eqn:least_squares}
\end{align}
where
\begin{align*}
	r =&\; 
	\begin{bmatrix}
		r[t_{}^{*}] & r[t_{}^{*}+1] & \cdots & r[\mathrm{T}] 
	\end{bmatrix}_{}^{\top} \;, \\
	e_{k}^{} =&\; 
	\begin{bmatrix}
		e_{k}^{}[t_{}^{*}] & e_{k}^{}[t_{}^{*}+1] & \cdots & e_{k}^{}[\mathrm{T}] 
	\end{bmatrix}_{}^{\top} \;.
\end{align*}
The control sequence, which results in an output sequence $\{y[t]\}_{t=t_{}^{*}}^{\mathrm{T}}$ that perfectly tracks the reference trajectory $\{r[t]\}_{t=t_{}^{*}}^{\mathrm{T}}$, can be computed via solving the linear equation:
\begin{align}
	u_{}^{\star} = G_{}^{-1}(r - d) \;. 
	\label{eqn:optimal_solution}
\end{align}
Without loss of generality, one can assume that ${x_{0}^{} = 0}$, and, equivalently, $d = 0$. Hence,~\eqref{eqn:optimal_solution} can be rewritten as
\begin{align}
	u_{}^{\star} = G_{}^{-1} r \;.
	\label{eqn:optimal_solution_2}
\end{align}
As argued in~\cite{AOR:96}, the direct inversion of $G$ is not practical in general since it requires having the exact information of $G$.  Besides, instead of inverting the entire matrix $G$, it is sufficient to compute the pre-image of $r$ under $G$.

\subsection{Gradient-based iterative learning algorithm} 
There are various techniques in the literature to solve the unconstrained optimization problem~\eqref{eqn:least_squares} iteratively. The gradient-based iterative learning control algorithm has been received an increasing attention (see, e.g.,~\cite{OHD:09, AOR:96, ChO:13}) due to its simplicity and light-weight computations compared to higher-order techniques. This algorithm generates the control inputs to be used in the next iteration using the relation:
\begin{align*}
	u_{k+1}^{} = u_{k}^{} + \gamma G_{}^{\top}e_{k}^{} \;,
\end{align*}
where $\gamma>0$ is the learning gain. Using this update law, the error evolves as 
\begin{align*}
	e_{k+1}^{} = \big(I - \gamma GG_{}^{\top} \big) e_{k}^{}.
\end{align*}
Using the norm inequality, provided in~\cite{HoJ:13}, we have:
\begin{align*}
	\parallel e_{k+1}^{} \parallel \; = \; \parallel (I - \gamma GG_{}^{\top})e_{k}^{} \parallel \; \leq \; \parallel I -\gamma GG_{}^{\top} \parallel \parallel e_{k}^{} \parallel \;.
\end{align*}
For minimum phase systems, the smallest singular value of the matrix $G$ is nonzero and if one picks $0<\gamma \leq \nicefrac{2}{\rho(GG^\top)}$, then $\Vert I - \gamma GG_{}^{\top} \Vert < 1$ holds. Consequently, $\Vert e_{k}^{} \Vert$ converges to zero linearly as $k \rightarrow \infty$.

It is worth noting that, for non-minimum phase systems, the matrix $G$ has some singular values that are very close to zero; therefore, it might be significantly ill-conditioned, leading $\Vert I-\gamma GG^\top\Vert$ to become nearly one. Taking into account the typical rounding errors that exists in numerical solvers, it is safe to assume that the matrix $G$ has zero singular values in order to avoid convergence issues due to mis-estimation of the optimal learning gain parameter.

\subsection{Trajectory tracking problem with sparsity constraint}
Trading off the accuracy of trajectory tracking for the sparsity in control signals amounts to solve
\begin{equation}\label{eq:cardinality_problem}
\begin{aligned}
	\text{minimize} &\quad \frac{1}{2}\parallel r - G u_{}^{} \parallel_{}^{2}  \\
	\text{subject to} &\quad \parallel T u_{}^{} \parallel_{0}^{} \leq M \\
	& \quad u \in \mathcal{U},
\end{aligned}
\end{equation}
where $M \leq N$ with $M\in\mathbb{N}_{0}^{}$ and  $T\in \mathbb{R}^{N-1\times N}$ is the difference matrix
\begin{align*}
	T = 
	\begin{bmatrix}
		-1 & 1 & 0 & 0 & \cdots & 0 & 0 \\
		0 	& -1 & 1 & 0 & \cdots & 0 & 0 \\
		0 	& 0 & -1 & 1 & \cdots & 0 & 0 \\
		\vdots 	& \vdots & \vdots & \vdots &  & \vdots & \vdots \\ 
		0 & 0 & 0 & 0 & \cdots & -1 & 1
	\end{bmatrix}\;,
\end{align*}
and $\mathcal{U}$ is a  compact and convex set which represents the practical limits on the input signal. 
For example, limits on the magnitude of input signal can be modeled by either by a box constraint or an upper bound on $\ell_\infty$-norm of control input $u$. With the cardinality constraint in~\eqref{eq:cardinality_problem}, one limits the number of  changes in control input values compared to the initial value $u[0]$, thereby promoting sparsity in the frequency of applying control input . 

However, due to the cardinality constraint, the problem~\eqref{eq:cardinality_problem} is non-convex and difficult to solve. A common heuristic method in the literature relies on the $\ell_1$-regularized problem
\begin{align}\label{eq:sparse_control_problem}
\begin{array}{ll}
	\text{minimize} & \frac{1}{2}\parallel r - G u_{}^{} \parallel_{}^{2} + \lambda\parallel T u_{}^{} \parallel_{1}^{} \;,\\
	\text{subject to} &  u \in \mathcal{U},
\end{array}
\end{align}
where the second term is referred as \emph{total variation} of signal $u$ and the problem~\eqref{eq:sparse_control_problem} is often called \emph{total variation denoising} in signal processing literature \cite{Rudin:92}.

\section{Sparse Iterative Learning Control}\label{sec:SILC}
In this section, we develop a first-order method to solve the regularized control problem, proposed in~\eqref{eq:sparse_control_problem}, iteratively. Our technique is based on backward-forward splitting method~\cite{Passty:79}, which is applied to the composite problem:
\begin{equation}\label{eq:problem_sum_two_function}
	\mbox{minimize} \; F(u):=f(u)+g(u) 
\end{equation}
where $f:\mathbb{E}\rightarrow \mathbb{R}$ is a differentiable convex function with Lipschitz continuous gradient $L$ satisfying 
\begin{equation}
\nonumber
\begin{aligned}
\Vert \nabla f(x) - \nabla f(y)\Vert \leq L \Vert x - y\Vert \quad \forall x, y \in \mathbb{E},
\end{aligned}
\end{equation}
and $g:\mathbb{E}\rightarrow (-\infty, +\infty]$ a proper closed convex function. 
Given a scalar $t>0$ the \emph{proximal} map associated to $g$ is defined as
\begin{equation}
\begin{aligned}
\mbox{prox}_t (g)(x) := \underset{u}{\mbox{argmin}} \left\{ g(u)+\dfrac{1}{2t}\Vert u-x\Vert^2 \right\}.
\end{aligned}
\end{equation} 
An important property of the proximal map is that this point (for any proper closed convex function $g$) is the unique solution to the associated minimization problem, and as a consequence, one has \cite[Lemma 3.1]{Beck:09}:
\begin{equation}
\begin{aligned}
(I+t\partial g)^{-1}(x)= {\mbox{prox}_t}(g)(x), \; \forall x\in \mathbb{E}.
\end{aligned}
\end{equation} 
This result  can be used to find the following optimality condition for \eqref{eq:problem_sum_two_function}
\begin{equation}\nonumber
\begin{aligned}
0 & \in t \nabla f(x^\star) + t \partial g(x^\star)\\
x^\star &= (I+t\partial g)^{-1}(I-t\nabla f) (x^\star),  
\end{aligned}
\end{equation}
and then further developed to obtain a backward-forward splitting based method  to solve \eqref{eq:problem_sum_two_function}
\begin{equation}\label{eq:backward_forward}
\begin{aligned}
x_{k+1} =& \mbox{prox}_{\gamma} (g) (x_{k}-\gamma \nabla f(x_{k}))\\
 =& \underset{x}{\mbox{argmin}}\; \left\{ g(x) + \dfrac{1}{2\gamma} \Vert x - (x_{k}-\gamma \nabla f(x_{k}))\Vert^2\right\}.
\end{aligned}
\end{equation}
For instance, if $f= \Vert A x-b\Vert^2$ and $g=\Vert x \Vert_1$, then the famous Iterative-Shrinkage-Thresholding Algorithm (ISTA) is recovered; see e.g.,~\cite{Comm:05}. We use backward-forward splitting method to solve~\eqref{eq:sparse_control_problem}. In particular, let 
\begin{equation}\label{eq:f_g}
\begin{aligned}
f(u) := \dfrac{1}{2}\Vert G u-r\Vert^2, \quad g(u) := \lambda \Vert T u \Vert_1 +\mathcal{I}_\mathcal{U}(u),
\end{aligned}
\end{equation}
with $\mathcal{I}_\mathcal{U}$ denoting the indicator function on $\mathcal{U}$; i.e., $\mathcal{I}_{\mathcal{U}}(u) = 0$ if $u\in \mathcal{U}$ and $\mathcal{I}_{\mathcal{U}}(u) =\infty$ otherwise.
Applying the backward-forward splitting, the sparse ILC update rule is given by
\begin{equation}\label{eq:BF_total_variation}
\begin{aligned}
&u_{k+1}  = \mbox{prox}_{\lambda/\gamma}(g)(u_{k}+\gamma G^\top e_{k})\\
&\; =\underset{u} {\mbox{argmin}} \; \left\{ \Vert T u\Vert_1+ \mathcal{I}_{\mathcal{U}}(u) + \dfrac{1}{2\lambda\gamma}
 \Vert u - (u_{k} + \gamma G^\top e_{k})\Vert^2 \right\} 
\end{aligned}
\end{equation}
Unlike the ISTA algorithm with simple $\ell_{1}^{}$-norm regularization, the sparse ILC iterations ~\eqref{eq:BF_total_variation} involve a proximal map that does not admit a closed-form solution. To tackle this problem, we develop an iterative \emph{dual}-based approach for the proximal step. In particular, we are interested in solving
\begin{equation}\label{eq:denoising_tv}
\begin{aligned}
\underset{u\in\mathcal U}{\mbox{minimize}} \; \left\{ \lambda \Vert T u\Vert_1 + \dfrac{1}{2}
 \Vert u - b\Vert^2 \right\} , \;
\end{aligned}
\end{equation}
by using a first-order method. We have the following result:

\begin{lemma} \label{lem:lemma_1}
Denote $\mathcal{P}^{n-1}\subset \mathbb{R}^{n-1}$ as the $n-1$ dimensional real space bounded by unit infinity norm (i.e., $p\in \mathcal{P}^{n-1}$  then $\Vert p \Vert_\infty \leq 1$) and $\mathcal{L}\in\mathbb{R}^{n\times n-1}$ given as
\begin{align*}
	\mathcal{L} = 
	\begin{bmatrix}
		1 & 0 & 0 & \cdots & 0 \\
		-1 & 1 & 0 & \cdots & 0 \\
		0 & -1 & 1 & \cdots & 0 \\
		\vdots & \vdots & \vdots &  & \vdots \\
		0 & 0 & 0 & \cdots & 1 \\
		0 & 0 & 0 & \cdots & -1 
	\end{bmatrix} \;.
\end{align*}
Let $p\in \mathcal{P}^{n-1}$ be the optimal solution of 
\begin{equation}\label{eq:minimize_tv_prox}
\begin{aligned}
 \underset{p\in \mathcal{P}^{n-1}}{\mbox{minimize}}&\; h(p):=-\Vert \Pi_\mathcal{U}(b-\lambda \mathcal{L}p) - (b-\lambda \mathcal{L}p) \Vert^2 \\
 &\; + \Vert b-\lambda \mathcal{L}p \Vert^2 \;.
\end{aligned}
\end{equation}
Then, the optimal solution of~\eqref{eq:denoising_tv} is given by
\begin{equation}\label{eq:opt_u_tv_denoising}
\begin{aligned}
u  & = \Pi_\mathcal{U}(b - \lambda \mathcal{L}p) .
\end{aligned}
\end{equation}
\end{lemma}

Next, we present the smoothness properties of~\eqref{eq:minimize_tv_prox}.

\begin{lemma} \label{lem:lemma_2}
The cost function in~\eqref{eq:minimize_tv_prox} is continuously differentiable, and its gradient is given by
\begin{equation}
\begin{aligned}
\nabla h(p):= -2\lambda \mathcal{L}^\top \Pi_\mathcal{U}(b-\lambda \mathcal{L} p) \;.
\end{aligned}
\end{equation}
Moreover, its Lipschitz constant is bounded by
\begin{multline}
\Vert \nabla h(p) - \nabla h(p\prime) \Vert \leq 2\lambda^2 \Vert \mathcal{L}^\top \Vert^2 \Vert p - p\prime\Vert \\
= 2\lambda^2 \rho(\mathcal{L}^\top \mathcal{L})  \Vert p - p\prime\Vert, \; \forall p, p\prime \in \mathcal{P}^{n-1}.
\end{multline}
Moreover, it follows ${\rho(\mathcal L^\top \mathcal L)\leq 4}$.
\end{lemma}

We are  now ready to form an accelerated projected gradient-based method to solve~\eqref{eq:minimize_tv_prox} and~\eqref{eq:opt_u_tv_denoising}. Algorithm~\ref{alg:1} solves the problem by employing a Nesterov-like acceleration applied to the dual domain. The technique offers a better convergence rate $O(1/k^2)$ as opposed to a gradient-based technique that converges at rate $O(1/k)$; see e.g., \cite{Beck:09}. 

\begin{algorithm}[H]
\caption{\textbf{Accelerated Projected Gradient} }
\label{alg:1}
\begin{algorithmic}[1]
\State Let $(N,\lambda, b)$ be given as input. Set $q_1 = 0$. 
\For{$k = 1,\dots, N$ compute}
\State \begin{equation} \nonumber
\begin{aligned}
p_k &= \Pi_{\mathcal{P}^{n-1}} \left[ q_k + \dfrac{1}{\lambda \rho(\mathcal L^\top \mathcal L)} \mathcal L^\top  \Pi_\mathcal{U}(b-\lambda \mathcal L q_k)\right]  \\
t_{k+1} &= \dfrac{1+\sqrt{1+4t_k^2}}{2}\\
q_{k+1} &= p_k + \dfrac{t_k-1}{t_{k+1}}(p_k - p_{k-1}) 
\end{aligned}
\end{equation}
where $\Pi_{\mathcal{P}^{n-1}}(x)_i = \dfrac{x_i}{\max\{1, \vert x_i\vert\}}$ for $i = 1,\dots, n-1$.
\EndFor
\State Return $(x^\star, p^\star) = (\Pi_\mathcal U (b- \lambda \mathcal L p_N), p_N)$.
\end{algorithmic}
\label{alg_summary}
\end{algorithm}

One can implement the sparse iterative learning control updates $u_{k+1}$ in~\eqref{eq:BF_total_variation} by first taking a gradient step on $u_k$ and then computing the proximal step via Algorithm~\ref{alg:1}. Algorithm~\ref{alg:2} proposes the gradient-based S-ILC method. 

\begin{algorithm}[H]
\caption{\textbf{Gradient-based S-ILC}}
\label{alg:2}
\begin{algorithmic}[1]
\State Let $(N_1, N_2, \lambda, \mathcal U)$ be given as input. Set $u_0, e_0$ to vector 0, and $\gamma = 1/\rho(G^\top G)$.
\For{$k = 1,\dots, N_1$ }
\State Set $b_k =  u_{k-1} +\gamma G^\top e_{k-1}$.
\State Run \textbf{Algorithm~\ref{alg:1}} with $(N_2,  \gamma\lambda, b_k, \mathcal U)$ and obtain $u_k\in \mathcal U$.

\State Apply $u_k$  to the plant and receive $e_k$.
\EndFor
\State Return $u^\star = u_{N_1}$.
\end{algorithmic}
\label{alg_summary}
\end{algorithm}

Next lemma confirms that the gradient-based S-ILC results in a non-increasing sequence.

\begin{lemma}\label{lem:3}
Consider the sequence $\{u_k\}_{k\geq 0}^{}$ generated by S-ILC. The associated functional values 
\begin{equation}\nonumber
\begin{aligned}
F(u_k):= \dfrac{1}{2}\Vert G u_k -r\Vert^2 + \lambda \Vert T u_k\Vert
\end{aligned}
\end{equation} 
is non-increasing. That is, for all $k\geq 1$,
\begin{equation}\nonumber
F(u_{k+1}) \leq F(u_{k}).
\end{equation}
\end{lemma}
Moreover, from \cite[Theorem 3.1]{Beck:09}, it follows that
\begin{equation}\nonumber
\begin{aligned}
F(u_k) - F(u^\star) \leq \dfrac{\rho(G^\top G) \Vert u_0 - u^\star \Vert}{2k},\; u_0\in \mathcal{U};
\end{aligned}
\end{equation}
where $u^\star$ is the optimal control input while $k\geq 1$ is the number of outer-loop iterations in the gradient-based S-ILC.

To accelerate the convergence of Algorithm~\ref{alg:2}, Nesterov-like iterations can be applied to its outer-loop. The straight-forward application of the Nesterov's method leads to the following updates:
\begin{equation}\label{eq:SILC_nest}
\begin{aligned}
b_k &=  y_k -\gamma G^\top (Gy_{k}-r),\\
y_{k+1} &= u_k + \dfrac{t_k-1}{t_{k+1}}(u_k - u_{k-1}),
\end{aligned}
\end{equation}
where $u_k$ is the control input obtained from inner-loop Algorithm~\ref{alg:1}. However, this requires the access to the reference signal $r$, which is not practical in ILC application.  We rewrite these updates to find a feasible formulation for ILC. Let $\Delta e_k :=e_k - e_{k-1}$ and $y_{k+1}:=u_k + \tau_{k+1} \Delta u_k$ where $\tau_{k+1} := ({t_k-1})/{t_{k+1}}$ and $\Delta u_{k} = u_k - u_{k-1}$. Now, we rewrite the $b_k$-th update in~\eqref{eq:SILC_nest} as 
\begin{equation}\nonumber
\begin{aligned}
b_{k+1} &= y_{k+1} -\gamma G^\top (G y_{k+1}-r)\\
& = u_k + \tau_{k+1} \Delta u_k -\gamma G^\top(G u_k + \tau_{k+1} G \Delta u_k -r),\\
& = u_k +\tau_{k+1} \Delta u_k + \gamma G^\top( e_k + \tau_{k+1} \Delta e_k ),
\end{aligned}
\end{equation}
which relates the auxiliary variable $b_k$ to -- the readily available -- control input and  error signals.

 \begin{algorithm}[H]
\caption{\textbf{Accelerated S-ILC}} 
 \label{alg:3}
\begin{algorithmic}[1]
\State Let $(N_1, N_2, \lambda, \mathcal U)$ be given as input. Set $u_{-1},u_0, e_{-1}, e_0$ to vector $0$, $ t_0 = 0$, $t_1 =1$, and $\gamma=1/\rho(G^\top G)$.
\For{$k = 1,\dots, N_1$}
\State Set 
\begin{equation} \nonumber
\begin{aligned}
t_{k} &= \dfrac{1}{2}+\dfrac{1}{2}\sqrt{1+4t_{k-1}^2},\; 
\quad \tau_{k} = \dfrac{t_{k-1}-1}{t_k},\\
b_k &=  u_{k-1} + \tau_k \Delta u_{k-1} + \gamma G^\top (e_{k-1}+ \tau_k \Delta e_{k-1}).
\end{aligned}
\end{equation}
 
\State Run \textbf{Algorithm~\ref{alg:1}} with $(N_2, \gamma\lambda, b_k, \mathcal U)$ and obtain $u_k\in \mathcal U$. 
\State Apply $u_k$  to the plant and receive $e_k$.
 
\EndFor
\State Return $u^\star = u_{N_1}$.
\end{algorithmic}
\label{alg_summary}
\end{algorithm}

Algorithm~\ref{alg:3} presents the accelerated Nesterov-like iterates to solve S-ILC. From \cite[Theorem 4.4]{Beck:09a}, it yields
\begin{equation}\nonumber
\begin{aligned}
F(u_k) - F(u^\star) \leq \dfrac{2\rho(G^\top G) \Vert u_0 - u^\star \Vert}{(k+1)^2},\; u_0\in \mathcal{U};
\end{aligned}
\end{equation}
where $k \geq 1$ is the outer-loop counter of Algorithm~\ref{alg:3}.

\section{Numerical Example}\label{sec:numerical}
\begin{figure}[t]
    \centering
    \includegraphics[width=0.2\textwidth]{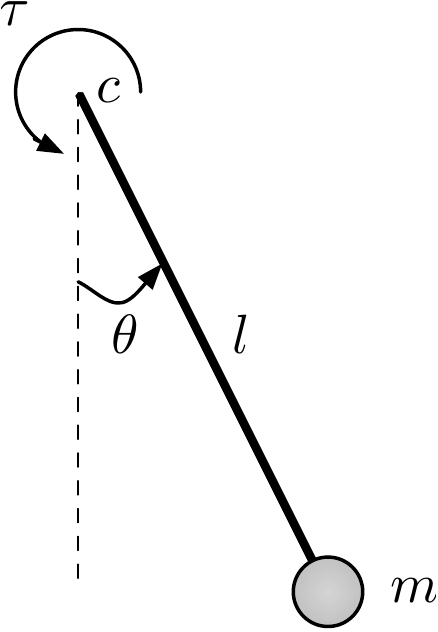} 
    \caption{A schematic drawing of the robot arm.}
    \label{fig:robot_arm}
\end{figure}

To demonstrate the effectiveness of S-ILC algorithms, we consider a robot arm ( see~\cite{VDS:11}) with one rotational degree-of-freedom as schematically shown in Fig.~\ref{fig:robot_arm}. The input is the torque $\tau$ applied to the arm at the joint and is limited to the range of $\pm \unit[12]{Nm}$, whereas the output $\theta$ is the angle of the arm measured as seen in Fig.~\ref{fig:robot_arm}. The dynamics of the robotic arm can be described by the following differential equation:
\begin{align}
	\ddot{\theta} = - \frac{g}{l}\sin\theta -  \frac{c}{ml_{}^{2}} \dot{\theta} + \frac{1}{ml_{}^{2}}\tau \;,
	\label{eqn:robot_arm_dynamic_model}
\end{align}
where the arm length is $l = \unit[1.0]{m}$, the payload mass is $m = \unit[1.0]{kg}$, the viscous friction coefficient is $c = \unitfrac[2.0]{Nms}{rad}$, and the gravitational acceleration is $g = \unitfrac[9.81]{m}{s^2}$. Changing the variables $x_{}^{(1)} \triangleq \theta$, $x_{}^{(2)} \triangleq \dot{\theta}$, $u \triangleq \tau$, and $y \triangleq \theta$, the nonlinear system~\eqref{eqn:robot_arm_dynamic_model} is sampled by using zero-order-hold and a sampling time of $T_{s}^{} = \unit[0.005]{s}$. The resulting discrete-time system becomes
\begin{align*}
	x_{}^{(1)}[t+1] =&\; x_{}^{(1)}[t] + T_{s}^{}x_{}^{(2)}[t] \;, \\
	x_{}^{(2)}[t+1] =&\; -\frac{g T_{s}^{}}{l}\sin(x_{}^{(1)}[t]) + \bigg( 1 - \frac{cT_{s}^{}}{ml_{}^{2}} \bigg)x_{}^{(2)}[t] \\
	&\; + \frac{T_{s}^{}}{ml_{}^{2}} u[t] \;, \\
	y[t] =&\; x_{}^{(1)}[t] \;.
\end{align*}

\begin{figure}[t]
    \centering
    \includegraphics[scale=1.00]{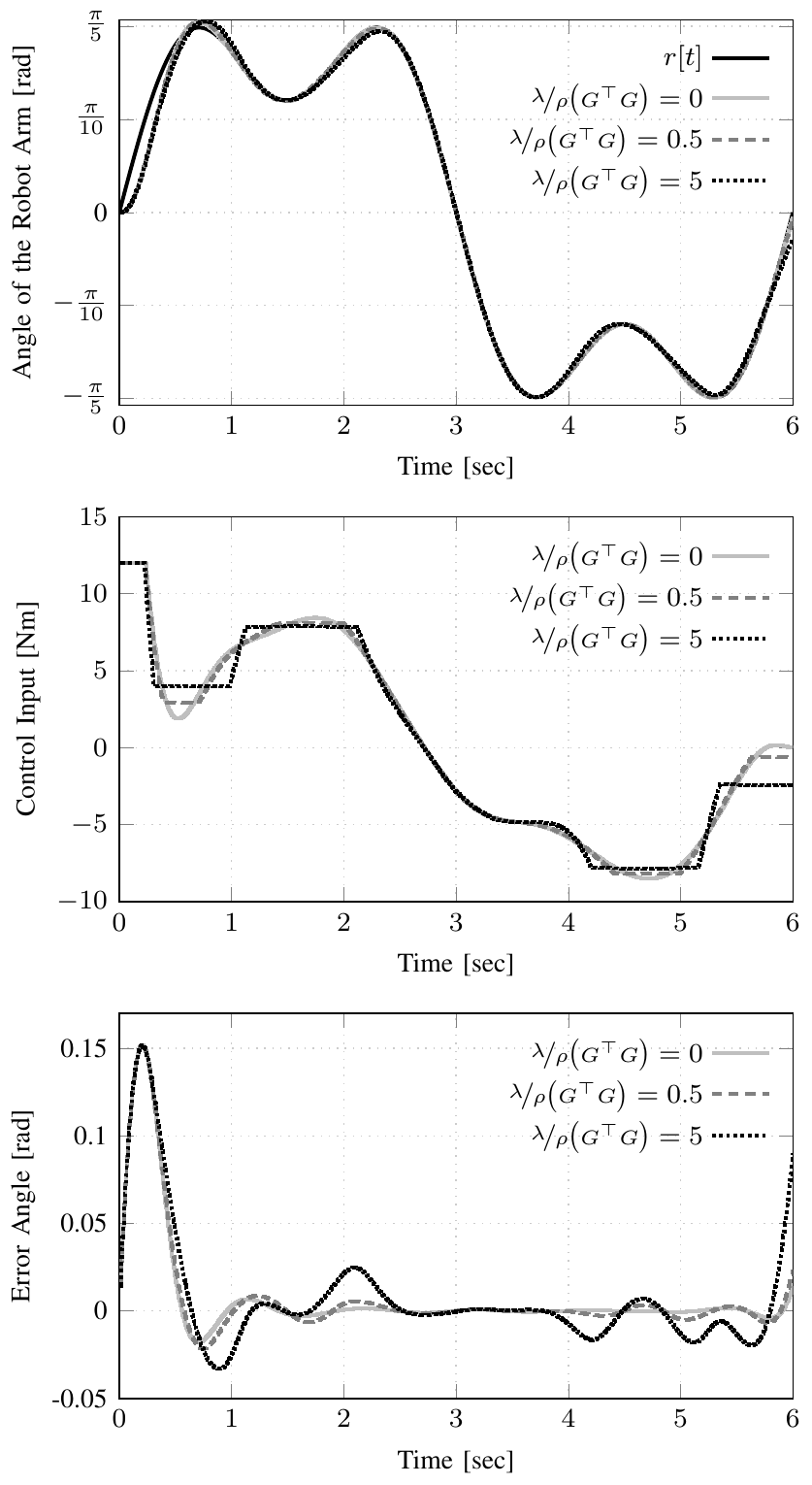} 
  	\caption{Tracking performance of the S-ILC algorithm for various values of the regularization parameter $\lambda$.}
  	\label{fig:results}
\end{figure}

To construct the gradient of the cost function~\eqref{eqn:least_squares}, which is equal to $G_{}^{\top}$, the discrete-time non-linear plant model is linearized around the stationary point $x_{}^{(1)} = \theta = 0$, resulting in the linear approximation:
\begin{align}
	x[t+1] &= 
	\begin{bmatrix}
		1 & T_{s}^{} \\ -\frac{g T_{s}^{}}{l} & 1 - \frac{cT_{s}^{}}{ml_{}^{2}} 
	\end{bmatrix}
	x[t] + 
	\begin{bmatrix}
		0 \\ \frac{T_{s}^{}}{ml_{}^{2}}	
	\end{bmatrix}
	u[t] \;, \\
	y[t] &= 
	\begin{bmatrix}
		1 & 0
	\end{bmatrix}
	x[t] \;.
\end{align}
Note that the linearized plant model is used to compute the gradient of the cost function~\eqref{eqn:least_squares}, whereas the nonlinear  model is employed in actual trials. The trial length is $\unit[6]{s}$ and the desired trajectory of the robot arm, illustrated in Fig.~\ref{fig:results}, is 
\begin{align*}
	r[t] = \frac{\pi}{5}\sin\bigg(\frac{\pi T_{s}^{} t}{3}\bigg) + \frac{2\pi}{25}\sin\big(\pi T_{s}^{} t\big) 
\end{align*}
for all $t \in \{0, 1, \cdots, 1200 \}$.

\begin{table}[ht]
\caption{Comparison of Regularization Parameters}
\centering
\begin{tabular}{c c c c}
\hline\hline
$\nicefrac{\lambda}{\rho\big(G_{}^{\top}G\big)}$ & $\parallel r - Gu \parallel_{2}^{}$ & $\parallel Tu \parallel_{1}^{}$ & $\parallel Tu \parallel_{0}^{}$  \\ [0.5ex] 
\hline
    $0$ 		& 1.0694 & 42.4495 & 1155 \\ 
    $0.5$ 	& 1.0845 & 38.0014 & 799 \\
    $2.5$ 	& 1.1406 & 34.5145 & 754 \\ 
    $5$ & 1.2117 & 33.0654 & 463 \\[1ex]
\hline
\end{tabular}
\label{tab:compare_parameters}
\end{table}

The simulation is carried out over 50 trials, and the results are displayed in Fig.~\ref{fig:results}. The optimization problem~\eqref{eq:sparse_control_problem} becomes a least square problem with a box constraint when $\lambda = 0$. The resulting control input sequence provides the smallest tracking error possible. As seen in Fig.~\ref{fig:results}, when the regularization parameter $\lambda$ increases, the control input sequence becomes more and more sparse at the expense of the increased tracking error. Similarly, Table~\ref{tab:compare_parameters} numerically illustrates the trade-off between the sparsity and the tracking performance. These experiments also demonstrate robustness against non-linearities of the plant. The change in the dynamics does not result in a divergence of the S-ILC algorithm.

\begin{figure}[t]
  \centering
  \includegraphics[scale=1.00]{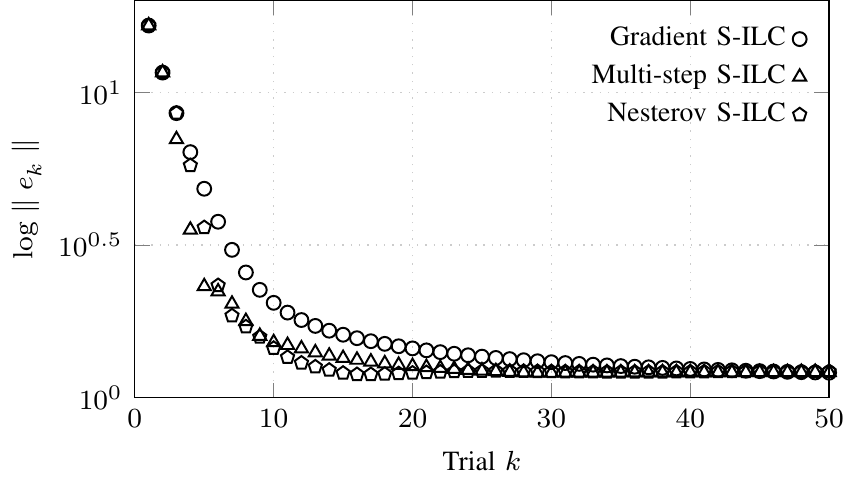} 
  \caption{The convergence of the error residual for different S-ILC methods.}
  \label{fig:error}
\end{figure}

Fig.~\ref{fig:error} shows the error decay rate of the gradient-based and accelerated S-ILC algorithms over $50$ trials. Moreover, we tried a multi-step technique called the heavy-ball  method which is obtained by adding a momentum term ${\beta (u_k - u_{k-1})}$ to the $b_{k+1}$-update in Algorithm~\ref{alg:2} where $\beta\in [0,1)$ is a scalar parameter. The superior convergence properties of the heavy-ball method compared to the gradient method is known for twice continuously differentiable cost functions~\cite{Polyak:87}. For the class of composite convex cost functions \eqref{eq:problem_sum_two_function}, however, the optimal algorithm parameters and associated convergence rate of the heavy-ball technique is still unknown~\cite{Ghadimi:15}. Here, we evaluated the heavy-ball algorithm with $\beta=0.4$. Numerical tests indicate that both Nesterov based (Algorithm~\ref{alg:3}) and heavy-ball methods improve the convergence of the gradient-based S-ILC Algorithm~\ref{alg:2}. 

It is noteworthy to mention that unlike the gradient-based algorithm, the accelerated and heavy-ball methods are not monotonic, that is, the function values~\eqref{eq:problem_sum_two_function} are not guaranteed to be non-increasing. In our evaluations, we observed if the inner-loop Algorithm~\ref{alg:1} is performed for a few iterations (so that it does not reach close to optimality), then the accelerated S-ILC problem (Algorithm~\ref{alg:3}) might become significantly non-monotonic to the point that it might diverge. Finding a monotone converging accelerated S-ILC algorithm is an interesting open problem.

\section{Conclusions}\label{sec:conclusion}
This paper has presented novel iterative learning algorithms to follow reference trajectories under a limited resource utilization. The proposed techniques promote sparsity by solving an $\ell_2$-norm optimization problem regularized by a total variation term. With proper selection of regularization parameter, our algorithms can strike a desirable trade-off between the accuracy of target tracking and the reduction of variations in actuation commands. Simulation results validated the efficacy of the proposed methods.

\section{Appendix}

\textit{Proof of Lemma~\ref{lem:lemma_1}:}
The result can be derived in a similar way as~\cite[Proposition 4.1]{Beck:09}.  Here, for completeness, we include the proof.
First, note  ${\vert x\vert =\underset{p} {\mbox{maximize}} \; \{px: \vert p\vert \leq 1\}}$,
and, similarly, $\Vert T u\Vert_1  = \sum_{i=1}^{n-1}\vert u_i - u_{i+1}\vert$ can be written as
\begin{equation}\nonumber
\begin{aligned} 
 \underset{p}{\mbox{maximize}}\;\{ \sum_{i=1}^{n-1}p_i (u_{i}-u_{i+1}): \vert p_i\vert \leq 1\} = \underset{p\in \mathcal P^{n-1}}{\mbox{maximize}}\; \mathcal L p^\top u.
\end{aligned}
\end{equation}
Accordingly, the problem~\eqref{eq:denoising_tv} becomes
\begin{equation}\nonumber 
\begin{aligned}
\underset{u\in \mathcal{U}}{\mbox{minimize}}\; \underset{p\in \mathcal{P}^{n-1}}{\mbox{maximize}}\; \dfrac{1}{2} \Vert u-b\Vert^2 + \lambda \mathcal{L}p^\top u.
\end{aligned}\end{equation}
Since this problem is convex in $u$ and concave in $p$, the order of minimization and maximization can be changed to obtain
\begin{equation}\nonumber 
\begin{aligned}
 \underset{p\in \mathcal{P}^{n-1}}{\mbox{maximize}}\; \underset{u\in \mathcal{U}}{\mbox{minimize}}\;\dfrac{1}{2} \Vert u-b\Vert^2 + \lambda \mathcal{L}p^\top u.
\end{aligned}\end{equation}
Using the basic relation 
 $$\Vert x-b \Vert^2 + 2c^\top x = \Vert x- b+c\Vert^2- \Vert b-c \Vert^2 + \Vert b\Vert^2,$$
 with $x = u$ and $c = \lambda \mathcal{L}p$ results in the equivalent form
\begin{equation}\nonumber 
\begin{aligned}
 \underset{p\in \mathcal{P}^{n-1}}{\mbox{maximize}}\; \underset{u\in \mathcal{U}}{\mbox{minimize}}\; \Vert u-(b-\lambda \mathcal{L}p)\Vert^2 - \Vert b- \lambda \mathcal{L}p\Vert^2 + \Vert b\Vert^2.
\end{aligned}\end{equation}
The optimal solution of the minimization problem, readily, is given by~\eqref{eq:opt_u_tv_denoising}. Instituting the optimal value of $u$, we arrive at the following dual problem
\begin{equation}
\begin{aligned}
 \underset{p\in \mathcal{P}^{n-1}}{\mbox{maximize}}&\; \Vert \Pi_\mathcal{U}(b-\lambda \mathcal{L}p) - (b-\lambda \mathcal{L}p) \Vert^2 - \Vert b-\lambda \mathcal{L}p \Vert^2 \;,
\end{aligned}
\end{equation}
which completes the proof.
\hfill $\blacksquare$

\textit{Proof of Lemma~\ref{lem:lemma_2}:}
Denote $s(x):= \dfrac{1}{2}\Vert x - \Pi_\mathcal{U} (x) \Vert^2$ and note that according to the \emph{proximal} map the following identity holds:
\begin{equation}\nonumber
\begin{aligned}
s(x) = \mbox{inf}_y \{ \Pi_\mathcal{U}(y)+\dfrac{1}{2}\Vert y - x\Vert^2 \}.
\end{aligned}
\end{equation}
From~\cite[Lemma 3.1]{Beck:09} it follows that $s(\cdot)$ is continuously differentiable with
\begin{equation}
\begin{aligned}
\nabla s (x): = x - \Pi_\mathcal{U}(x).
\end{aligned}
\end{equation}
The gradient of $h(p)$ then reads
\begin{equation}\nonumber
\begin{aligned}
\nabla h(p)& = \nabla ( -2s(b-\lambda \mathcal{L}p)+ \Vert b-\lambda \mathcal{L}p \Vert^2)\\
& =  2\lambda \mathcal{L}^\top (b-\lambda \mathcal{L}p - \Pi_\mathcal{U}(b-\lambda \mathcal{L}p)) - 2\lambda \mathcal{L}^\top (b-\lambda \mathcal{L}p)\\
& = -2\lambda \mathcal{L}^\top \Pi_\mathcal{U}(b-\lambda \mathcal{L} p).
\end{aligned}
\end{equation}
For any $p, p\prime \in \mathcal{P}^{n-1}$ we have
\begin{equation}\nonumber
\begin{aligned}
\Vert \nabla h(p) - \nabla &h(p\prime)\Vert \\
& =  \Vert 2\lambda \mathcal{L}^\top (\Pi_\mathcal{U}(b-\lambda \mathcal{L}p) - \Pi_\mathcal{U}(b-\lambda \mathcal{L}p\prime)) \Vert\\
& \leq 2\lambda \Vert \mathcal L^\top\Vert \Vert \Pi_\mathcal{U}(b-\lambda \mathcal{L}p) - \Pi_\mathcal{U}(b-\lambda \mathcal{L}p\prime)\Vert\\
& \leq 2\lambda^2 \Vert\mathcal L^\top\Vert  \Vert\mathcal L (p-p\prime) \Vert\\
& \leq 2\lambda^2 \Vert \mathcal{L}^\top\Vert \Vert \mathcal L\Vert \Vert p-p\prime\Vert \\ 
& = 2\lambda^2 \rho(\mathcal{L}^\top \mathcal{L}) \Vert p -p\prime\Vert.
 \end{aligned}
\end{equation}
The matrix $\mathcal L^\top \mathcal L \in \mathbb{R}^{n-1\times n-1}$ is given by
\begin{align*}
	\mathcal{L}_{}^{\intercal}\mathcal{L} = 
	\begin{bmatrix}
     2   &  -1  &   0  &     \cdots   &  0   &  0 \\
    -1   &   2   & -1  &      \cdots   &  0   &  0 \\ 
     0   &  -1   &  2  &     \cdots  &   0 &    0 \\ 
     \vdots  &   \vdots   &  \vdots   &     & \vdots   &  \vdots \\
     0 &    0  &   0 &     \cdots    &  2  &  -1 \\
     0  &   0   &  0  &   \cdots    & -1 &    2
	\end{bmatrix} \;.
\end{align*}
Following the Gershgorin circle theorem \cite{HoJ:13} one concludes the largest eigenvalue of this matrix follows $\rho(\mathcal L^\top \mathcal L) \leq \max_i\; \sum_{j=1}^{n-1} {\vert [\mathcal L^\top \mathcal L]_{ij} \vert }\leq 4.$ 
\hfill $\blacksquare$

\textit{Proof of Lemma~\ref{lem:3}:}
Let $f(u)$ and $g(u)$ be defined as \eqref{eq:f_g} and 
\begin{equation}\nonumber
P_{\lambda/\rho}(u):=\mbox{prox}_{\lambda/\rho(G^\top G)} (g) (u - \dfrac{1}{\rho(G^\top G)}G^\top(G u-r)).
\end{equation} 
Then the S-ILC algorithm (with  converging inner-loop and no modeling error) can be rewritten as 
\begin{equation}\nonumber
\begin{aligned}
u_{k+1} = P_{\lambda/\rho} (u_k).
\end{aligned}
\end{equation} 
For a convex Lipschitz continuous gradient  function $f$ and convex function $g$, define \cite{Beck:09a}:
 $$Q_L(x,y) = f(y)+\langle \nabla f(y), x-y\rangle + \dfrac{L}{2}\Vert x-y\Vert^2 + g(x)$$
Then, it can be seen that $P_{\lambda/\rho} (u) = \underset{x}{\mbox{argmin}} \, Q_{\lambda/\rho}(x, u)$. Furthermore, we have
\begin{equation}\nonumber
\begin{aligned}
F(u_k)&\geq Q_{\lambda/\rho}(u_k, u_k)  \geq Q_{\lambda/\rho}(P_{\lambda/\rho}(u_k), u_k)\\
& =f(u_k) + \langle \nabla f(u_k), P_{\lambda/\rho}(u_k)-u_k\rangle \\
&+ \dfrac{\rho(G^\top G)}{2}\Vert P_{\lambda/\rho}(u_k)-u_k\Vert^2+ g(P_{\lambda/\rho}(u_k)) \\
& \geq f(P_{\lambda\rho}(u_k)) + g(P_{\lambda/\rho}(u_k)) \\
&= F(u_{k+1}),
\end{aligned}
\end{equation}
where the last inequality holds for Lipschitz continuous $f$.

\hfill $\blacksquare$

\balance
\bibliographystyle{styles/IEEEtran}
\bibliography{bibILC}

\end{document}